\documentclass[11pt]{article}

\usepackage{amsmath,amssymb,amsthm}
\usepackage{amsfonts}
\usepackage[mathscr]{eucal}
\begin{document}
\newtheorem{thm}{Theorem}
\numberwithin{thm}{section}
\newtheorem{lemma}[thm]{Lemma}
\newtheorem{remark}{Remark}
\newtheorem{corr}[thm]{Corollary}
\newtheorem{proposition}{Proposition}
\newtheorem{theorem}{Theorem}[section]
\newtheorem{deff}[thm]{Definition}
\newtheorem{case}[thm]{Case}
\newtheorem{prop}[thm]{Proposition}
\numberwithin{equation}{section}
\numberwithin{remark}{section}
\numberwithin{proposition}{section}
\newtheorem{corollary}{Corollary}[section]
\newtheorem{others}{Theorem}
\newtheorem{conjecture}{Conjecture}\newtheorem{definition}{Definition}[section]
\newtheorem{cl}{Claim}
\newtheorem{cor}{Corollary}
\newcommand{\ds}{\displaystyle}

\newcommand{\stk}[2]{\stackrel{#1}{#2}}
\newcommand{\dwn}[1]{{\scriptstyle #1}\downarrow}
\newcommand{\upa}[1]{{\scriptstyle #1}\uparrow}
\newcommand{\nea}[1]{{\scriptstyle #1}\nearrow}
\newcommand{\sea}[1]{\searrow {\scriptstyle #1}}
\newcommand{\csti}[3]{(#1+1) (#2)^{1/ (#1+1)} (#1)^{- #1
 / (#1+1)} (#3)^{ #1 / (#1 +1)}}
\newcommand{\RR}[1]{\mathbb{#1}}
\thispagestyle{empty}
\begin{titlepage}
\title{\bf Laws of the iterated logarithm for $\alpha$-time Brownian motion}
\author{Erkan Nane\thanks{ Supported in part by NSF Grant \# 9700585-DMS }\\
Department of Mathematics\\
Purdue University\\
West Lafayette, IN 47906 \\
enane@math.purdue.edu}
\maketitle
\begin{abstract}
\noindent {\it  We introduce a class of iterated processes called
$\alpha$-time Brownian motion for $0<\alpha \leq 2$. These are
obtained by taking Brownian motion and replacing the time
parameter with a symmetric $\alpha$-stable process. We prove a
Chung-type law of the iterated logarithm (LIL) for these processes
which is a generalization of LIL proved  in \cite{hu} for iterated
Brownian motion. When $\alpha =1$ it takes the following form
$$
\liminf_{T\to\infty}T^{-1/2}(\log \log T) \sup_{0\leq t\leq
T}|Z_{t}|=\pi^{2}\sqrt{\lambda_{1}} \ \ \ a.s.
$$
where $\lambda_{1}$ is the first eigenvalue for the Cauchy process
 in the interval $[-1,1].$
We also define the local time $L^{*}(x,t)$ and range
$R^{*}(t)=|\{x:\ Z(s)=x \text{ for some } s\leq t\}|$ for these
processes for $1<\alpha <2$. We prove that there are universal
constants $c_{R},c_{L}\in (0,\infty) $ such that
$$
\limsup_{t\to\infty}\frac{R^{*}(t)}{(t/\log \log
t)^{1/2\alpha}\log \log t}= c_{R} \ \ a.s.
$$
$$
 \liminf_{t\to\infty} \frac{\sup_{x\in \RR{R}}L^{*}(x,t)}{(t/\log
\log t)^{1-1/2\alpha}}= c_{L} \ \ a.s.
$$}
\end{abstract}
\end{titlepage}
\section{Introduction}

In recent years, several iterated processes received much research
interest from many mathematicians,  see \cite{allouba1, burdzy1,
chuan, deblassie, nane1, nane2, shi-yor, xiao} and references
there in. Inspired by these results, we introduce a new class of
iterated processes called $\alpha$-time Brownian motion for
$0<\alpha \leq 2$. These are obtained by taking Brownian motion
and replacing the time parameter with  a symmetric $\alpha$-stable
process. For $\alpha =2$, this is the iterated Brownian motion of
Burdzy \cite{burdzy1}. One of the main differences of these
iterated processes and Brownian motion is that they  are not
Markov or Gaussian. However,  for $\alpha =1,2$ these processes
have connections with partial differential operators as described
in \cite{allouba1, nane3}.

To define $\alpha$-time Brownian motion, let $X_{t}$ be a
two-sided Brownian motion on $\RR{R}$. That is, $\{X_{t}:\ t\geq
0\}$ and $\{X_{-t}:\ t\leq 0\}$ are two independent copies of
Brownian motion starting from $0$.  Let $Y_{t}$ be a real-valued
symmetric $\alpha$-stable process, $0<\alpha \leq 2$,  starting
from $0$ and independent of $X_t$. Then $\alpha$-time Brownian
motion $Z_{t}$ is defined by
\begin{equation}\label{definition}
Z_{t}\equiv X(Y_{t}), \ \ \ t\geq 0.
\end{equation}
It is easy to verify that $Z_{t}$ has stationary increments and is
a self-similar process of index $1/2\alpha$.  That is, for every
$k>0$, $\{Z_{t}:\ t\geq 0 \}$ and $\{k^{-1/ 2\alpha}Z_{kt}:\ t\geq
0\}$ have the same finite-dimensional distributions. We refer to
Taqqu \cite{taqqu} for relations of self-similar stable processes
to physical quantities.  The $\alpha$-time Brownian motion is an
 example of nonstable self-similar processes.

%
Our aim in this paper is two-fold. Firstly, we will be interested
in the path properties of the process defined in
(\ref{definition}). Since this process is not Markov or Gaussian,
it is of interest to see how the lack of independence of
increments affect the asymptotic behavior. Secondly, we will
define the local time $L^{*}(x,t)$ for this process for $1< \alpha
<2$. We will prove the joint continuity  of the local time and
extend LIL of Kesten \cite{kesten} to these processes. We also
obtain an LIL for the range of these processes.

In the first part of the paper, we will be interested in proving a
``liminf" law of the iterated logarithm of the Chung-type for
$Z_{t}$. The study of this type of LIL's was initiated by Chung
\cite{chung} for Brownian motion $W_{t}$. He proved that
$$
\liminf_{T\to\infty}( T^{-1} \log \log T)^{1/2} \sup_{0\leq t\leq
T}|W_{t}|=\frac{\pi}{8^{1/2}} \ \ a.s.
$$

This LIL was extended to several other processes later
including symmetric $\alpha$-stable processes  $Y_{t}$ by Taylor \cite{taylor} in the following form
$$
\liminf_{T\to\infty}(T^{-1}\log \log T)^{1/\alpha}\sup_{0\leq
t\leq T}|Y_{t}|=(\lambda_{\alpha})^{1/\alpha}\ \ a.s.
$$
where $\lambda_{\alpha}$ is the first eigenvalue of the fractional
Laplacian $(-\Delta )^{\alpha /2}$ in $[-1,1]$.

One then wonders if a Chung-type LIL holds for the composition of
symmetric stable processes. Although these processes are not Markov or Gaussian, this has been
achieved for the composition of two Brownian motions, the so called
iterated Brownian motion, which is the case of $\alpha =2$ proved by  Hu,
Pierre-Loti-Viaud, and Shi in
\cite{hu}. They showed that
\begin{equation}\label{IBM-chung}
\liminf_{T\to\infty}T^{-1/4}(\log \log T)^{3/4} \sup_{0\leq t\leq
T}|S_{t}^{1}| =\left( \frac{3\pi^{2}}{8}\right)^{3/4} \ \
a.s.
\end{equation}
with $S_{t}^{1}=X(W_{t})$ denoting iterated Brownian motion, where
$X$ is a two-sided Brownian motion and $W$ is another Brownian
motion independent of $X$. This is the definition of iterated
Brownian motion used by Burdzy \cite{burdzy1}.

Inspired by the above mentioned extensions of the Chung's LIL
we extend the above results to  composition of a Brownian motion and a
symmetric $\alpha$-stable process.

\begin{theorem}\label{maintheorem}
Let $\alpha\in (0,2]$ and let $Z_{t}$ be  the $\alpha$-time Brownian
motion as defined in (\ref{definition}). Then we have
\begin{equation}\label{mainlimit}
\liminf_{T\to\infty}T^{-1/2\alpha}(\log \log
T)^{(1+\alpha)/(2\alpha)} \sup_{0\leq t\leq T}|Z_{t}|=D_{\alpha}\
\ a.s.
\end{equation}
where $D_{\alpha}=C_{\alpha}^{(1+\alpha)/2\alpha}$,
$C_{\alpha}=(\pi^{2}/8)^{\alpha/(1+\alpha)}
(1+\alpha)(2^{\alpha}\lambda_{\alpha})^{1/(1+\alpha)}(\alpha)^{-\alpha/(1+\alpha)}$.
\end{theorem}

A Chung-type LIL  has also been established for other versions of iterated
 Brownian motion (see \cite{chuan}, \cite{koslew}) as follows:

\begin{equation}\label{IBM-chung2}
\liminf_{T\to\infty}T^{-1/4}(\log \log T)^{3/4} \sup_{0\leq t\leq
T}|S_{t}| =\frac{3^{3/4}\pi^{3/2}}{2^{11/4}} \ \ a.s.
\end{equation}
with $S_{t}\equiv W(|\hat{W_{t}}|)$ denoting another version of
iterated Brownian motion, where $W$ and $\hat{W}$ are independent
real-valued standard Brownian motions, each starting from $0$. For
a generalization of this result to $\alpha$-time Brownian motions
we define the  process
\begin{equation}\label{definition1}
Z^{1}_{t}\equiv X(|Y_{t}|), \ \ \ t\geq 0.
\end{equation}
for Brownian motion $X_{t}$ and symmetric
$\alpha$-stable process  $Y_{t}$ independent of $X$, each starting
from $0$, $0<\alpha \leq 2$. For this process we have

\begin{theorem}\label{maintheorem1}
Let $\alpha\in (0,2]$ and let $Z^{1}_{t}$ be the $\alpha$-time Brownian
motion as defined in (\ref{definition1}). Then we have
\begin{equation}\label{mainlimit1}
\liminf_{T\to\infty}T^{-1/2\alpha}(\log \log
T)^{(1+\alpha)/(2\alpha)} \sup_{0\leq t\leq
T}|Z^{1}_{t}|=D^{1}_{\alpha}\ \ a.s.
\end{equation}
where $D^{1}_{\alpha}=(C^{1}_{\alpha})^{(1+\alpha)/2\alpha}$,
$C^{1}_{\alpha}=(\pi^{2}/8)^{\alpha/(1+\alpha)}
(1+\alpha)(\lambda_{\alpha})^{1/(1+\alpha)}(\alpha)^{-\alpha/(1+\alpha)}$.
\end{theorem}

We note that the constants appearing  in (\ref{mainlimit}) and
(\ref{mainlimit1}) are different. The main reason for this is that
the process $Z_{t}$ have three independent processes $\{X_{t}:\
t\geq 0\}$,  $\{X_{-t}:\ t\leq 0\}$ and $Y$, while the process
$Z^{1}_{t}$ does not have a contribution from $\{X_{-t}:\ t\leq
0\}$. The proof of Theorem \ref{maintheorem1} follows the same
line of proof of Theorem \ref{maintheorem}, except for the small
deviation probability estimates for $Z^{1}_{t}$ we use Theorem
\ref{small-ball1}.

The motivation for the study of Local times of $\alpha$-time
Brownian motion came from the results of  Cs\'{a}ki,
Cs\"{o}rg\"{o}, F\"{o}ldes, and R\'{e}v\'{e}sz \cite{CCFR} and Shi
and Yor \cite{shi-yor} about Kesten--type  laws of iterated
logarithm for iterated Brownian motion. The study of this type of
LIL's was initiated by Kesten \cite{kesten}. Let $B_{t}$ be a
Brownian motion, $L(x,t)$ its local time at $x$. Then Kesten
showed
\begin{equation}\label{kest-law1}
\limsup_{t\to\infty}\frac{L(0,t)}{\sqrt{2t\log \log t}}=
\limsup_{t\to\infty}\frac{\sup_{x\in \RR{R}}L(x,t)}{\sqrt{2t\log
\log t}}=1 \ \ a.s.
\end{equation}
and
\begin{equation}\label{kest-law2}
\liminf_{t\to\infty}\frac{\sup_{x\in \RR{R}}L(x,t)}{\sqrt{t /\log
\log t}}=c \ \  a.s. \ \ 0<c<\infty.
\end{equation}

These types of laws  were generalized later to symmetric stable processes of index $\alpha \in (1,2)$.
More specifically, Donsker and
Varadhan \cite{don-var} generalized  (\ref{kest-law1})
and Griffin \cite{griffin} generalized (\ref{kest-law2}) to symmetric stable processes.


More recently, Kesten--type LIL's were extended to iterated
Brownian motion $S$. Let $L_{S}(x,t)$ be the local time of
$S_{t}=W_{1}(|W_{2}(t)|)$, with $W_{1}$ and $W_{2} $ independent
standard real-valued Brownian motions. (\ref{kest-law1}) was
extended to the IBM case by Cs\'{a}ki, Cs\"{o}rg\"{o}, F\"{o}ldes,
and R\'{e}v\'{e}sz \cite{CCFR}  and Xiao \cite{xiao}.  This result
asserts that there exist (finite) universal constants $c_{1}>0$
and $c_{2}>0$ such that

\begin{equation}\label{ibm-kest1}
c_{1}\leq \limsup_{t\to\infty} \frac{L_{S}(0,t)}{t^{3/4}(\log \log
t)^{3/4}} \leq c_{2} \ \ \ a.s.
\end{equation}
(\ref{kest-law2}) was extended to IBM case by Cs\'{a}ki,
Cs\"{o}rg\"{o}, F\"{o}ldes, and R\'{e}v\'{e}sz \cite{CCFR} and Shi
and Yor \cite{shi-yor}. This result asserts that there exist
universal constants $c_{3}>0$ and $c_{4}>0$ such that
\begin{equation}\label{ibm-kest2}
c_{3}\leq \liminf_{t\to\infty} t^{-3/4}(\log \log t)^{3/4}
\sup_{x\in \RR{R}} L_{S}(x,t) \leq c_{4} \ \ a.s.
\end{equation}

Inspired by the definition of local time of IBM in \cite{CCFR}, we
define the local time of $\alpha$-time Brownian motion defined in
(\ref{definition1}) as follows:
\begin{eqnarray}
L^{*}(x,t)& = & \int_{0}^{\infty} \bar{L}_{2}(u,t)d_{u}L_{1}(x,u)\nonumber\\
& = & \int_{0}^{\infty} (L_{2}(u,t)+
L_{2}(-u,t))d_{u}L_{1}(x,u)\label{alpha-local},
\end{eqnarray}
  where $L_{1}$,  $L_{2}$  and $\bar{L}_{2}$ denote, respectively, the local times of $X$, $Y$
and $|Y|$. A similar definition can be given for $\alpha$-time Brownian
motion defined in (\ref{definition}) using the ideas in
\cite{burdzy-khos}.

In \S3, we will extend (\ref{ibm-kest2}) to $\alpha$-time
Brownian motion. We will also obtain partial results towards
extension of (\ref{kest-law1}). However, our results does not
imply an extension of LIL in (\ref{kest-law1}) and are far from
optimal yet, and leave many problems open.
 These results will follow from the study of L\'{e}vy classes for the local
time of $Z$ and $Z^{1}$. We extend (\ref{ibm-kest2}) as follows.
\begin{theorem}\label{local-liminf}
There exists a universal constant $c_{L} \in (0,\infty)$ such that
$$
 \liminf_{t\to\infty} \frac{\sup_{x\in
\RR{R}}L^{*}(x,t)}{(t/\log \log t)^{1-1/2\alpha}}= c_{L}\ \ a.s.
$$
\end{theorem}
A similar result holds also for the local time of the process
defined in (\ref{definition}). The usual LIL or Kolmogorov's LIL
for Brownian motion which replaces the time parameter was used
essentially in the results in \cite{CCFR} and \cite{shi-yor} to
prove Kesten's LIL for iterated Brownian motion. However, there
does not exist an LIL of this type for symmetric $\alpha$-stable
process which replaces the time parameter in the definition of
$\alpha$-time Brownian motion. To overcome this difficulty we show
in Lemma \ref{distr-measure} that the LIL for the range process of
symmetric $\alpha$-stable process suffices to prove Theorem
\ref{local-liminf}.

We also obtain usual LIL for the range of $\alpha$-time Brownian
motion. Then, we use it with a particular case of occupation times
formula to obtain Kesten's LIL for these processes. This
  is also essential
in the study of some of L\'{e}vy classes of local time of $Z^{1}$.

\begin{theorem}\label{range-upper}
There exists a universal constant $c_{R} \in (0,\infty)$ such that
$$
\limsup_{t\to\infty}\frac{R^{*}(t)}{(t/\log \log
t)^{1/2\alpha}\log \log t}= c_{R} \ \ a.s.
$$
where $R^{*}(t)=|\{ x:\ Z^{1}(s)=x \ \text{for some }\ s\leq
t\}|$.
\end{theorem}
A similar result holds also for range of the process defined in
(\ref{definition}).

Our proofs of Theorems \ref{maintheorem} and \ref{maintheorem1} in
this paper follow the proofs  in \cite{hu} making necessary
changes
 at crucial points. In studying the local time of
 $\alpha$-time Brownian motion we use the ideas we learned from
 Cs\'{a}ki, Cs\"{o}rg\"{o}, F\"{o}ldes,
and R\'{e}v\'{e}sz \cite{CCFR} and Shi and Yor \cite{shi-yor} with
necessary changes in the use of usual LIL of the range processes.
Our proofs differ from theirs since there does not exist usual LIL
for symmetric $\alpha$-stable process. To overcome this difficulty
we show that the usual LIL for the range of symmetric
$\alpha$-stable process suffice for our results, see Lemma
\ref{distr-measure}. We also adapt the arguments of Griffin
\cite{griffin} to our case to prove the usual LIL for the range of
$\alpha$-time Brownian motion. Our proofs differs from his in that
$\alpha$-time Brownian motion does not have independent
increments. So we disjointify the range of symmetric
$\alpha$-stable process to get independent increments, see Lemma
\ref{prob-bound}. The paper is organized as follows. We give the
proof of Theorem \ref{maintheorem} in \S2. The local time and the
range of
 $\alpha$-time Brownian motion are studied in \S3.


\section{Chung's LIL for $\alpha$-time Brownian motion}

We will prove Theorem \ref{maintheorem} in this section. Section
2.1 is devoted to the preliminary lemmas about the small deviation
probabilities. In section 2.2 we prove the lower bound in Theorem
\ref{maintheorem}. Upper bound is proved in section 2.3.
\subsection{Preliminaries}

In this section we give some definitions and preliminary lemmas which will be used
in the proof of the main result.

A real-valued symmetric stable process $Y_{t}$ with index
$\alpha\in (0,2]$ is the process with stationary independent
increments whose transition density
$$
p_{t}^{\alpha}(x,y)=p^{\alpha}(t,x-y), \ \ \ \ \ (t,x,y)\in
(0,\infty)\times \RR{R}^{n}\times \RR{R}^{n},
$$
is characterized by the Fourier transform
$$
\int_{\RR{R}^{n}}
e^{iy.\xi}p^{\alpha}(t,y)dy=\exp(-t|\xi|^{\alpha}), \ \ \ \ \ t>0,
\xi \in \RR{R}^{n}.
$$
The process has right continuous paths, it is rotation and
translation invariant.

The following lemma gives the small ball probabilities for the process $\sup_{0\leq
t\leq 1}|Y_{t}|$.
\begin{lemma}[Mogul'skii, 1974, \cite{mogul}]\label{M}
Let $0<\alpha \leq 2$ and let $Y_{t}$ be a symmetric
$\alpha$-stable process. Then
$$
\lim_{\epsilon\to 0^{+}}\epsilon^{\alpha}\log P[\ \sup_{0\leq
t\leq 1}|Y_{t}|\leq \epsilon \ ]=-\lambda_{\alpha},
$$
where $\lambda_{\alpha}$ is the first eigenvalue of the fractional
Laplacian operator in the interval $[-1,1]$.
\end{lemma}

This is an equivalent statement of the fact that
$$
\lim_{t\to\infty} P[\tau>t]=-\lambda_{\alpha},
$$
due to scaling property of $\sup_{0\leq
t\leq 1}|Y_{t}|$, where
$
\tau=\inf \{s:\ |Y_{s}|\geq 1\}
$
is the first exit time of the interval $[-1,1]$

Let $R=\sup_{0\leq t\leq 1}Y_{t}\ - \ \inf_{0\leq t\leq 1}Y_{t} $ be the range of $Y_{t}$.
The following is a special case of Theorem 2.1 in \cite{Chenli}.

\begin{theorem}[Mogul'skii, 1974, \cite{mogul}]\label{range}
$$\lim_{\epsilon\to 0^{+}}\epsilon^{\alpha}\log P[\ R \leq \epsilon \ ]=-2^{\alpha}\lambda_{\alpha}.$$
\end{theorem}

We use the following theorem  (Kasahara \cite[Theorem 3]{kasahara}
and Bingham, Goldie and Teugels \cite[p. 254]{bgt}) to find the  asymptotics of the Laplace transform of $R$ below.

\begin{theorem}[de Bruijn's Tauberian Theorem]
Let $X$ be a positive random variable such that for some positive
$B_{1},\ B_{2}$ and $p$,
$$-B_{1}\leq \liminf_{x\rightarrow 0} x^{p}\log P[X \leq x] \leq
\limsup_{x\rightarrow 0} x^{p}\log P[X \leq x]\leq -B_{2}.
$$
Then
\begin{eqnarray}
\  & \ &  -(p+1)(B_{1})^{1/(p+1)}p^{-p/(p+1)}
\leq \liminf_{\lambda \rightarrow \infty} \lambda^{-p/(p+1)}\log E e^{- \lambda X}\nonumber \\
\ & \ & \leq  \limsup_{\lambda \rightarrow \infty}
\lambda^{-p/(p+1)}\log E e^{- \lambda X} \leq
-(p+1)(B_{2})^{1/(p+1)}p^{-p/(p+1)} .\nonumber
\end{eqnarray}

\end{theorem}
From de Bruijn's Tauberian theorem and Theorem \ref{range} we have
\begin{lemma}\label{bruijn-range}
$$\lim_{\lambda\to\infty}\lambda^{-\alpha/(1+\alpha)}\log
E[e^{-\lambda
R}]=-(1+\alpha)(2^{\alpha}\lambda_{\alpha})^{1/(1+\alpha)}\alpha^{-\alpha
/(1+\alpha)}.$$
\end{lemma}

The following theorem gives the small ball deviation probabilities for the process $Z_{t}$ defined in (\ref{definition}).
\begin{theorem}\label{small-ball}
We have
$$
\lim_{u\to 0}u^{2\alpha /(1+\alpha)}\log P[\ \sup_{0\leq t\leq 1}
|Z_{t}|\leq u \ ]=-(\pi^{2}/8)^{\alpha /(1+\alpha)}A_{\alpha},
$$
where
$A_{\alpha}=(1+\alpha)(2^{\alpha}\lambda_{\alpha})^{1/(1+\alpha)}\alpha^{-\alpha
/(1+\alpha)}$.
\end{theorem}
\begin{proof}
Let $X_{t}$ be a Brownian motion. From a well-known formula (see
Chung \cite{chung}):
$$
P[\ \sup_{0\leq t \leq 1} |X_{t}|\leq u \
]=\frac{4}{\pi}\sum_{k=1}^{\infty}\frac{(-1)^{k-1}}{2k-1}\exp
\left[ -\frac{(2k-1)^{2}\pi^{2}}{8u^{2}}\right],
$$
we get that, for all $u>0$,
\begin{equation}\label{twoside}
\frac{2}{\pi}\exp\left( -\frac{\pi^{2}}{8u^{2}}\right)\leq P[\
\sup_{0\leq t \leq 1} |X_{t}|\leq u \ ]\leq
\frac{4}{\pi}\exp\left( -\frac{\pi^{2}}{8u^{2}} \right).
\end{equation}
Let $Z_{t}=X(Y_{t})$ be the $\alpha$-time Brownian motion and let
\begin{equation}
S(t)\equiv \sup_{0\leq s \leq t} Y_{s}, \ \ I(t)\equiv \inf_{0\leq
s \leq t} Y_{s},
\end{equation}
then

$
P[\ \sup_{0\leq t\leq 1} |Z_{t}|\leq u \ ]
$
\begin{eqnarray}
& =& P[\ \sup_{0\leq t\leq S(1)} |X_{t}|\leq u, \ \sup_{I(1)\leq
t\leq 0} |X_{t}|\leq u \ ] \nonumber\\
&=&E\left[ P\left( \sup_{0\leq t\leq 1} |X_{t}|\leq
\frac{u}{\sqrt{S(1)}}|Y \right) P\left( \sup_{0\leq t\leq 1}
|X_{t}|\leq \frac{u}{\sqrt{-I(1)}}|Y \right)\right]\nonumber\\
&\leq  &
\frac{16}{\pi^{2}}E\exp\left(-\frac{\pi^{2}(S(1)-I(1))}{8u^{2}}\right).
\end{eqnarray}
This last inequality follows from the second part of
(\ref{twoside}). Similarly the first part of (\ref{twoside}) gives
us a lower bound, with $4\pi^{-2}$ instead of $16\pi^{-2}$. Now
the proof follows from the given inequalities and Lemma
\ref{bruijn-range}
\end{proof}
The following theorem gives the small ball deviation probabilities for the process $Z^{1}_{t}$ defined in (\ref{definition1}).
\begin{theorem}\label{small-ball1}
We have
$$
\lim_{u\to 0}u^{2\alpha /(1+\alpha)}\log P[\ \sup_{0\leq t\leq 1}
|Z^{1}_{t}|\leq u \ ]=-(\pi^{2}/8)^{\alpha /(1+\alpha)}A^{1}_{\alpha},
$$
where
$A^{1}_{\alpha}=(1+\alpha)(\lambda_{\alpha})^{1/(1+\alpha)}\alpha^{-\alpha
/(1+\alpha)}$.
\end{theorem}
\begin{proof}
Let $M(t)\equiv \sup_{0\leq s\leq t}|Y_{s}|$.
The proof follows the same line of the proof of Theorem \ref{small-ball}, except at the end we have

$
P[\ \sup_{0\leq t\leq 1} |Z^{1}_{t}|\leq u \ ]
$
\begin{eqnarray}
& =& P[\ \sup_{0\leq t\leq M(1)} |X_{t}|\leq u,  ] \nonumber\\
&=&E\left[ P\left( \sup_{0\leq t\leq 1} |X_{t}|\leq
\frac{u}{\sqrt{M(1)}}|Y \right) \right]\nonumber\\
&\leq  &
\frac{4}{\pi}E\exp\left(-\frac{\pi^{2}M(1)}{8u^{2}}\right),
\end{eqnarray}
and similarly a lower bound with $2/\pi$ instead of $4/\pi$. Then we use Lemma \ref{M} and de Bruijn's Tauberian theorem.
\end{proof}

\begin{lemma}\label{two-sided-sup}
For all $0<a\leq b,$ $u>0$ we have
$$
P[\ a< \sup_{0\leq t \leq u}|Z_{t}| <b \ ]\leq (b/a -1)^{2}.
$$
\end{lemma}
\begin{proof}
The proof follows from the proof of Lemma 4.1 in \cite{hu}.
\end{proof}
The following proposition is the combination of two propositions
in \cite{bertoin}, which are stated as Proposition 2 on page 219
and Proposition 4 on page 221.
\begin{proposition}\label{sup-y}
Let $Y_{t}$ be a symmetric $\alpha$-stable process. Let
$$ S(t)=\sup_{0\leq s\leq t} Y_{s}.$$
 Then there exists $k_{1},
k_{2}>0$ such that
$$
\lim_{x\to\infty} x^{\alpha}P[Y_{1}>x]=\lim_{x\to\infty}
x^{\alpha} P[S(1)>x]= k_{1},
$$
and
$$
\lim_{x\to 0^{+}}x^{-\alpha /2}P[S(1)<x]= k_{2} .
$$
\end{proposition}

We will use following versions of Borel-Cantelli lemmas in our
proofs.
\begin{lemma}[Borel-Cantelli Lemma 1]
Let $E_{1}, E_{2},\cdots$ be a sequence of events (sets) for which
$\sum_{n=1}^{\infty}P[E_{n}]<\infty$. Then
$$
P[E_{n} \ i.o]=P[\cap_{n=1}^{\infty}\cup_{i=n}^{\infty}E_{i}]=0,
$$
i.e with probability $1$ only a finite number of events $E_{n}$
occur simultaneously.
\end{lemma}

Since the process $Z_{t}$ does not have independent increments we
have to use another version of the Borel-Cantelli lemma which is
due to Spitzer \cite{spitzer}.
\begin{lemma}[Borel-Cantelli Lemma 2, p. 28 in \cite{revez}]\label{liminf-prob}
Let $E_{1}, E_{2},\cdots$ be a sequence of events (sets) for which
$$
\sum_{n=1}^{\infty}P[E_{n}]=\infty
 \ \ \text{and} \ \liminf_{n\to\infty}\sum_{i=1}^{n}\sum_{j=1}^{n}P[E_{i}E_{j}]\left/
\left( \sum_{i=1}^{n}P[E_{i}]\right)^{2} \right. \leq c \ \ (c\geq
1).
$$
Then  $$ P[\cap_{n=1}^{\infty}\cup_{k=n}^{\infty}E_{k}]=P[E_{n}\
i.o ]\geq 1/c.$$
\end{lemma}

\subsection{Proof of the lower bound}

The lower bound is easier as always. We use Theorem \ref{small-ball}.
Let
$$
C_{\alpha}=(\pi^{2}/8)^{\alpha/(1+\alpha)}
(1+\alpha)(2^{\alpha}\lambda_{\alpha})^{1/(1+\alpha)}(\alpha)^{-\alpha/(1+\alpha)},
$$
be the small deviation probability limit for $\sup_{0\leq t \leq
T}|Z_{t}|$ given in Theorem \ref{small-ball}.
 For every fixed $\epsilon >0$, it follows from
Theorem \ref{small-ball} that, for $T$ sufficiently large, we have
\newline

       $ P\left[ T^{-1/2\alpha} (\log \log T)^{(1+\alpha)/2\alpha}\sup_{0\leq
t \leq T}|Z_{t}|\leq
(1-\epsilon)^{(1+\alpha)/\alpha}C_{\alpha}^{(1+\alpha)/2\alpha}\right]
$
\begin{eqnarray}
&=& P\left[ \sup_{0\leq t \leq 1}|Z_{t}|\leq
(1-\epsilon)^{(1+\alpha)/\alpha}C_{\alpha}^{(1+\alpha)/2\alpha}(\log
\log T)^{-(1+\alpha)/2\alpha} \right ]\nonumber\\
&\leq & \exp\left[ -(1-\epsilon) (1-\epsilon)^{-2}C_{\alpha}
C_{\alpha}^{-1}\log \log T \right]\nonumber\\
& = & \exp\left[ -\frac{1}{(1-\epsilon)} \log \log T
\right].\nonumber
\end{eqnarray}
Taking a fixed rational number $a>1$ and $T_{k}=a^{k}$ gives that
\begin{eqnarray}
& & \sum_{k\geq 1}P[ T_{k}^{-1/2\alpha} (\log \log
T_{k})^{(1+\alpha)/2\alpha} \sup_{0\leq t \leq T_{k}}|Z_{t}|\nonumber  \\
& & \ \ \ \ \ \ \ \ \ \ \ \  \leq
(1-\epsilon)^{(1+\alpha)/\alpha}C_{\alpha}^{(1+\alpha)/2\alpha}]<+\infty.\nonumber
\end{eqnarray}
 It follows from Borel-Cantelli lemma, by letting $\epsilon\to 0$,
 that
 \begin{equation}\label{lower-for-k}
\liminf_{k\to + \infty}T_{k}^{-1/2\alpha} (\log \log
T_{k})^{(1+\alpha)/2\alpha}\sup_{0\leq t \leq T_{k}}|Z_{t}|\geq
C_{\alpha}^{(1+\alpha)/2\alpha}= D_{\alpha}\ a.s.
 \end{equation}
 Since for every $T>0$, there exists $k\geq 0$ such that $T_{k}\leq T<
 T_{k+1}$, we have
 $$ T^{-1/2\alpha} (\log \log T)^{(1+\alpha)/2\alpha}\sup_{0\leq t
\leq T}|Z_{t}| \geq  T_{k+1}^{-1/2\alpha} (\log \log
T_{k})^{(1+\alpha)/2\alpha}\sup_{0\leq t \leq T_{k}}|Z_{t}|$$
$$\ \ \ \ \  \ \ \ \ \ \ \ \ \ \ \ \ \ \ \ \ \ \ \ \ \ \ \ \ \ \ \ \
 =
a^{-1/2\alpha}T_{k}^{-1/2\alpha} (\log \log
 T_{k})^{(1+\alpha)/2\alpha}\sup_{0\leq t \leq
T_{k}}|Z_{t}|,$$
which together with (\ref{lower-for-k}), yields the lower bound,
as the rational number $a>1$ can be arbitrarily close to $1$.

\subsection{Proof of the upper bound}

We follow the steps in the
proof of Lemma 4.2 in \cite{hu}. Let $\epsilon
>0$ be fixed. For notational simplicity, we use the following in
the sequel
\begin{eqnarray}
& T_{k}& = \exp(k\log k) \nonumber\\
&a_{k}&= (1+3\epsilon)^{(1+\alpha)/\alpha}C_{\alpha}^{(1+\alpha)/2\alpha}
T_{k}^{1/2\alpha}(\log \log T_{k})^{-(1+\alpha)/2\alpha}
\nonumber\\
&B_{k}&=\{ \sup_{0\leq t\leq T_{k}}|Z_{t}|\leq a_{k}\}.  \nonumber
\end{eqnarray}
It follows from Theorem \ref{small-ball} that there exists
$k_{o}(\epsilon)$, depending only on $\epsilon$, such that for
every $k> k_{o}(\epsilon)$, we have
\begin{eqnarray}
P(B_{k}) & \geq & \exp \left( - (1+3\epsilon)
C_{\alpha}(1+3\epsilon )^{-2}C_{\alpha}^{-1} \log \log
T_{k}\right)\nonumber\\
&=& \exp\left(  -\frac{1}{1+3\epsilon}\log \log
T_{k}\right)\nonumber\\
& \geq & k^{-1/(1+2\epsilon)},\nonumber
\end{eqnarray}
which yields existence of positive constants $C=C(\epsilon)$ and
$N=N(\epsilon)$ such that for every $n>N$,
\begin{equation}\label{greater}
\sum_{k=1}^{n}P(B_{k})\geq Cn^{2\epsilon /(1+2\epsilon)}.
\end{equation}
We now establish the following
\begin{lemma}\label{liminf}
We have
\begin{equation}\label{gen-liminf}
\liminf_{n\to\infty}\sum_{i=1}^{n}\sum_{j=1}^{n}P[B_{i}B_{j}]\left/
\left( \sum_{i=1}^{n}P[B_{i}]\right)^{2} \right. \leq 1.
\end{equation}
\end{lemma}

\begin{proof}[Proof of Lemma \ref{liminf}]
Let $K>0$ be a constant such that
$$K\geq 1/\alpha(3(1+2\epsilon)/\epsilon -2),$$
 and let
$n_{\epsilon}=[n^{\epsilon  /(1+2\epsilon)}]$ (with $[x]$ denoting
the integer part of $x$). We set furthermore
$$
E_{1}=\{ (i,j): \ 1\leq i, j \leq n, \ |i-j|\leq 19 \}
$$
$$
 E_{2}=\{
(i,j): \ n_{\epsilon} \leq i, j \leq n, \ |i-j|\geq 20 \}.
$$
Since by (\ref{greater}),
$$
\sum_{i=1}^{n_{\epsilon}}\sum_{j=1}^{n}P[B_{i}B_{j}]\left/ \left(
\sum_{i=1}^{n}P[B_{i}]\right)^{2} \right. \leq n_{\epsilon}\left/
\sum_{i=1}^{n}P[B_{i}] \right. \leq C^{-1} n^{-\epsilon
/(1+2\epsilon)},
$$
and
$$
\sum \sum_{(i,j)\in E_{1}} P[B_{i}B_{j}]\left/ \left(
\sum_{i=1}^{n}P[B_{i}]\right)^{2} \right. \leq 38\left/
\sum_{i=1}^{n}P[B_{i}] \right.\leq 38 C^{-1}n^{-2\epsilon
/(1+2\epsilon)},
$$
it suffices to prove that
\begin{equation}\label{E2-liminf}
\liminf_{n\to\infty} \sum \sum_{(i,j)\in E_{2}}
P[B_{i}B_{j}]\left/ \left( \sum_{i=1}^{n}P[B_{i}]\right)^{2}
\right. \leq 1.
\end{equation}
Let
$$
S(t) =\sup_{0\leq s\leq t}Y_{s},\ \ \
I(t)=\inf_{0\leq s\leq t}Y_{s},
$$
$$
F(x)=P[\sup_{0\leq t\leq 1}|X_{t}|\leq x],\ \ \ \ \forall x>0.
$$
Then for all $i<j$ and all positive numbers $ p_{1}<p_{2}, \
q_{1}<q_{2}$,
\begin{eqnarray}
& & P[B_{i}B_{j}|S(T_{i})=p_{1}, S(T_{j})=p_{2}, I(T_{i})=q_{1}, I(T_{j})=q_{2}]\nonumber\\
 & & =P[\sup_{-q_{1}\leq t\leq p_{1}}|X_{t}|\leq a_{i},\ \sup_{p_{1}\leq t\leq p_{2}}|X_{t}|\leq a_{j},
 \  \sup_{-q_{2}\leq t\leq -q_{1}}|X_{t}|\leq a_{j} ]\nonumber\\
 & &=P[\sup_{0\leq t\leq p_{1}}|X_{t}|\leq a_{i},\ \sup_{p_{1}\leq t\leq p_{2}}|X_{t}|\leq
 a_{j}]\nonumber\\
 & &\ \ \ \times P[\sup_{0\leq t\leq q_{1}}|X_{t}|\leq a_{i},\ \sup_{q_{1}\leq t\leq q_{2}}|X_{t}|\leq
 a_{j}].\nonumber
\end{eqnarray}
Notice that for all $x>0$ and $y>0$, \newline
 $ P[\sup_{0\leq
t\leq p_{1}}|X_{t}|\leq x,\ \sup_{p_{1}\leq t\leq
p_{2}}|X_{t}|\leq
 y]$
 \begin{eqnarray}
& & \leq P[\sup_{0\leq t\leq p_{1}}|X_{t}|\leq x]\sup_{|u|\leq
x}P[\sup_{0\leq t\leq p_{2}-p_{1}}|X_{t}+u|\leq
 y]\label{u-added}\\
 & & =P[\sup_{0\leq t\leq p_{1}}|X_{t}|\leq x]P[\sup_{0\leq t\leq p_{2}-p_{1}}|X_{t}|\leq
 y]\label{u-added-result}\\
 & & =F(xp_{1}^{-1/2})F(y(p_{2}-p_{1})^{-1/2}).\nonumber
 \end{eqnarray}
 The equation (\ref{u-added}) is from the Markov property of
 Wiener processes, and equation (\ref{u-added-result}) is due to a
 general property of Gaussian measures (see, e.g. Ledoux and Talagrand \cite[p.73]{ledoux})

 It follows that
 \begin{eqnarray}
 P[B_{i}]&=&
 E[F(a_{i}S^{-1/2}(T_{i}))F(a_{i}(-I(T_{i}))^{-1/2})],\nonumber\\
 P[B_{i}B_{j}]& \leq & E\left[  F\left( \frac{a_{i}}{\sqrt{S(T_{i})}}\right)
 F\left( \frac{a_{i}}{\sqrt{-I(T_{i})}}\right)\right. \nonumber\\
& &\ \ \ \left.  \times F\left(
\frac{a_{j}}{\sqrt{S(T_{j})-S(T_{i})}}\right)
 F\left( \frac{a_{j}}{\sqrt{I(T_{i})-I(T_{j})}}\right) \right].\nonumber
 \end{eqnarray}
 Let $\mathcal{F}_{t}\equiv \{Y_{s}, s\leq t \}$ and $M(t)\equiv \sup_{0\leq s\leq t}|Y_{s}|
 $. Let $f(t)=\exp(Kt)$.  By noticing that $ S(T_{j})-S(T_{i})$
 (resp. $I(T_{i})-I(T_{j})$) is bounded below by the positive part of
 $\sigma^{+}$ (resp. $\tau^{+}$) of
 \begin{eqnarray}
 \sigma & \equiv  & \sup_{T_{i}\leq t \leq
 T_{j}}(Y_{t}-Y_{T_{i}})-2M(T_{i})\nonumber\\
 & &\ \ (resp.\ \  \tau \equiv -\inf_{T_{i}\leq t \leq
 T_{j}}(Y_{t}-Y_{T_{i}})-2M(T_{i})),
 \end{eqnarray}
\begin{eqnarray}
(i.e.\ \sup_{T_{i}\leq t \leq
 T_{j}}(Y_{t}-Y_{T_{i}})-2M(T_{i})^{+}& \leq & (\sup_{T_{i}\leq t \leq
 T_{j}}(Y_{t})+|Y_{T_{i}}|-2M(T_{i}))^{+}\nonumber\\
 & \leq & M(T_{j})-M(T_{i}) ),\nonumber
\end{eqnarray}
 we obtain that
\begin{eqnarray}
& & E[F( a_{j}(S(T_{j})-S(T_{i}))^{-1/2})
 F( a_{j}(I(T_{j})-I(T_{i}))^{-1/2})|\mathcal{F}_{T_{i}}]\nonumber\\
& & \ \ \ \ \leq E[F( a_{j}(\sigma^{+})^{-1/2})
 F( a_{j}(\tau^{+})^{-1/2})|\mathcal{F}_{T_{i}}]\nonumber\\
 & &\ \ \ \ \leq \sup _{0\leq x\leq 2M(T_{i})}E[F(
 a_{j}((S(T_{j}-T_{i})-x)^{+})^{-1/2})\label{equiv-distr}\\
 & &\ \ \ \  \ \ \times F( a_{j}((I(T_{j}-T_{i})-x)^{+})^{-1/2})]\nonumber\\
& & \ \ \ \ \leq  1_{\{M(T_{i})>T_{i}^{1/\alpha}f(\log \log T_{i})\}} +
E[F( a_{j}(\mu^{+})^{-1/2})
 F( a_{j}(\nu^{+})^{-1/2})],\nonumber
\end{eqnarray}
with
$$\mu =S(T_{j}-T_{i})-2 T_{i}^{1/\alpha}f(\log \log T_{i})\
and \  \nu = -I(T_{j}-T_{i})-2 T_{i}^{1/\alpha}f(\log \log
T_{i}).$$ In the inequality (\ref{equiv-distr}) we use the fact
that $\alpha$-stable process has  stationary independent
increments.

As $F$ is always between $0$ and $1$, we get that
\begin{eqnarray}
P[B_{i}B_{j}]&\leq & P[ M(T_{i})>T_{i}^{1/\alpha}f(\log \log
T_{i})]\nonumber\\
& & \ \ + P[B_{i}]E[F( a_{j}(\mu^{+})^{-1/2})
 F( a_{j}(\nu^{+})^{-1/2})]\nonumber\\
 &=& P[ M(T_{i})>T_{i}^{1/\alpha}f(\log \log
T_{i})]\nonumber \\
& & \ \ +P[B_{i}]E\left[F\left(
\frac{a_{j}}{\sqrt{(S(1)-\theta)^{+}(T_{j}-T_{i})^{1/\alpha}}}\right)
\right.\nonumber\\
 & & \ \ \ \left. \times F\left(
\frac{a_{j}}{\sqrt{(-I(1)-\theta)^{+}(T_{j}-T_{i})^{1/\alpha}}}\right)\right].\label{theta-inequality}
\end{eqnarray}
The identity (\ref{theta-inequality}) is due to the scaling
property of $\alpha $-stable process, with
$$
\theta = 2(T_{j}-T_{i})^{-1/\alpha}T_{i}^{1/\alpha}f(\log \log
T_{i}).
$$
Now by using Proposition \ref{sup-y},  $\lim_{x\to\infty}
x^{\alpha}P[S(1)>x]= k_{1}$, we get that, if $i\geq
n_{\epsilon}=[n^{\epsilon  /(1+2\epsilon)}] $, then since $ f(\log
\log T_{i})$ is large
\begin{eqnarray}
P[ M(T_{i})>T_{i}^{1/\alpha}f(\log \log T_{i})]& \leq & 2P[
S(1)>f(\log \log T_{i})]\nonumber\\
& \leq  & 4k_{1} (f(\log \log T_{i}))^{-\alpha}\nonumber\\
&=& 4 k_{1}(i \log i)^{-\alpha K}\nonumber\\
& \leq & 4k_{1}n^{-\alpha K \epsilon /(1+2\epsilon )},
\end{eqnarray}
i.e.

\begin{eqnarray}
& &\sum_{i=n_{\epsilon}}^{n} \sum_{j=1}^{n}  P[
M(T_{i}) > T_{i}^{1/\alpha}f(\log \log T_{i})]\left/  \left(
\sum_{i=1}^{n}P[B_{i}] \right)^{2}\right.\nonumber\\
& & \ \ \ \ \leq  4k_{1}C^{-2}n^{2-(\alpha K+2)\epsilon/(1+2\epsilon)} \nonumber\\
& &\ \ \ \  \leq  4k_{1}C^{-2}n^{-1}, \label{lessthan1}
\end{eqnarray}
as $K\geq 1/\alpha(3(1+2\epsilon)/\epsilon -2)$. On the other
hand, for $j\geq i+20$
$$
\theta (  T_{j}/T_{i})^{1/2\alpha}\leq 2(i\log
i)^{K}/(j^{(j-i)/2}-1)^{1/\alpha}\leq 2C_{0}j^{-(j-i)/5\alpha}\leq
2C_{0}j^{-20/5\alpha},
$$
which is small for the range of $j$ we consider (if needed we can take $j\geq i+20 +C(K)$, where $C(K)$ is a constant multiple of $K$). Since from
Proposition \ref{sup-y}, for $x$ close to $0$,
$$
P[S(1)<x]\leq (1+\epsilon)k_{2}x^{\alpha /2},
$$
we have
\begin{eqnarray}
P[S(1)-\theta < (1-(T_{j}/T_{i})^{1/2\alpha})S(1)]&=&
P[S(1)<\theta (  T_{j}/T_{i})^{1/2\alpha} ]\nonumber\\
&\leq  & 2k_{2}(\theta (  T_{j}/T_{i})^{1/2\alpha})^{\alpha
/2}\nonumber\\
& \leq  & 2k_{2}( 2C_{0}j^{-20/5\alpha})^{\alpha /2}\nonumber\\
&\leq &C_{1} j^{-2},\nonumber
\end{eqnarray}
with some universal constant $C_{1}$. This inequality holds for
$-I(1)$ instead of $S(1)$ as well, since symmetric $\alpha$-stable
process is symmetric. Therefore for all $(i,j)\in E_{2}$, $C_{2}$
being a universal constant, we have
\begin{eqnarray}
& & E\left[F\left(
  \frac{a_{j}}{\sqrt{(S(1)-\theta)^{+}(T_{j}-T_{i})^{1/\alpha}}}\right)
F\left(
\frac{a_{j}}{\sqrt{(-I(1)-\theta)^{+}(T_{j}-T_{i})^{1/\alpha}}}\right)\right]\nonumber\\
& & \ \ \leq 2C_{1}j^{-2} +E\left[ F\left(
\frac{a_{j}}{\sqrt{(S(1)G_{ij}}}\right) F\left(
\frac{a_{j}}{\sqrt{(-I(1)G_{ij}}}\right)\right]\nonumber\\
& &\ \ \ \ \ \ \ \ (G_{ij}\equiv
(T_{j}-T_{i})^{1/\alpha}(1-(T_{i}/T_{j})^{1/2\alpha}))\nonumber\\
& & \ \  \leq 2C_{1}j^{-2} + E\left[ F\left(
\frac{a_{j}(1+C_{2}j^{-(j-i)/5\alpha})}{\sqrt{S(1)T_{j}^{1/\alpha}}}\right)F\left(
\frac{a_{j}(1+C_{2}j^{-(j-i)/5\alpha})}{\sqrt{-I(1)T_{j}^{1/\alpha}}}\right)\right]\nonumber\\
& &\ \ \ \ \ \ \ \  ( since\ \ \sqrt{T_{j}^{1/\alpha} / G_{ij}}\leq
1+C_{2}j^{-(j-i)/5\alpha} )\nonumber \\
& & \ \ \leq 2C_{1}j^{-2}+ P[\sup_{0 \leq t \leq T_{j}} |Z_{t}|\leq
a_{j}(1+C_{2}j^{-2})]\nonumber\\
& & \ \  \leq  2C_{1}j^{-2} +P[B_{j}]+P[a_{j}\leq \sup_{0\leq t \leq  T_{j}}
|Z_{t}|\leq a_{j}(1+C_{2}j^{-2}) ]\nonumber\\
& & \ \  \leq 2C_{1}j^{-2} +P[B_{j}]+C_{2}^{2}j^{-4},\nonumber
\end{eqnarray}
where the last inequality follows from Lemma \ref{two-sided-sup}.
Combining this with (\ref{theta-inequality}), (\ref{lessthan1})
and (\ref{greater}) gives that
\begin{eqnarray}
 & & \sum \sum_{(i,j)\in E_{2}} P[B_{i}B_{j}]\left/ \left(
\sum_{i=1}^{n}P[B_{i}]\right)^{2} \right.\nonumber\\
& & \ \ \ \ \leq 4k_{1}C^{-2}n^{-1}+1
+(2C_{1}+C_{2}^{2})\sum_{i=1}^{n}P[B_{i}]\sum_{j=1}^{n}j^{-2}\left/\left(
\sum_{i=1}^{n}P[B_{i}]\right)^{2}\right.\nonumber\\
& & \ \ \ \ \leq 1+
4k_{1}C^{-2}n^{-1}+\pi^{2}(2C_{1}+C_{2}^{2})(6C)^{-1}n^{-2\epsilon/(1+2\epsilon)},\nonumber
\end{eqnarray}
which yields (\ref{E2-liminf}).
\end{proof}
Since $\sum_{k=1}^{\infty}P[B_{k}]=\infty$, it follows from
(\ref{gen-liminf}) and a well-known version of Borel-Cantelli
lemma (Lemma \ref{liminf-prob} above) that
$P[\limsup_{k\to\infty}B_{k}]=1$ which implies the upper bound in
Theorem \ref{maintheorem}.

\section{Local time of $\alpha$-time Brownian motion }

In this section we give the definition of the local time of
$\alpha$-time Brownian motion and prove its joint continuity. In
section 3.0.1 we prove a lemma which is crucial in the proofs of
the main theorems. Sections 3.1-3.3  and section 3.5 give a study
of the L\'{e}vy classes for the local time. In section 3.4 we
prove an LIL for the range of $\alpha$-time Brownian motion.

Let $L_{1}(x,t)$ be the local time of Brownian motion, and
$L_{2}(x,t)$ be the local time of symmetric $\alpha$-stable
process for $1<\alpha \leq 2$ (see \cite{revez} for the properties
of the local time of Brownian motion and see \cite{griffin} and
references there in for the properties of the local time of
symmetric $\alpha$-stable processes). Let $f$, $x\in\RR{R}$, be a
locally integrable real-valued function. Then
\begin{equation}\label{stable-integral}
\int_{0}^{t}f(W_{i}(s))ds=\int_{-\infty}^{\infty}f(x)L_{i}(x,t)dx,
\ \ i=1,2,
\end{equation}
for $W_{1}$ a standard Brownian motion and $W_{2}$ a symmetric
stable process.

 Then we define the local time of the $\alpha$-time Brownian
motion as
\begin{equation}\label{local-definition}
L^{*}(x,t):=\int_{0}^{\infty}\bar{L}_{2}(s,t)d_{s}L_{1}(x,s),\ \
x\in \RR{R},\ \ t\geq 0,
\end{equation}
where $\bar{L}_{2}(x,t):=L_{2}(x,t)+L_{2}(-x,t) $, $x\geq 0$.

We prove next the joint continuity of $L^{*}(x,t)$ and establish
the occupation times formula for $Z^{1}_{t}$.
\begin{proposition}
There exists an almost surely  jointly continuous family of "local
times", $\{L^{*}(x,t): \ t\geq 0, \ x\in \RR{R} \}$, such that for
all Borel measurable integrable functions, $f: \RR{R}\to \RR{R}$
and all $t\geq 0$,
\begin{equation}\label{iterated-int}
\int_{0}^{t}f(Z^{1}(s))ds=\int_{0}^{t}f(X(|Y(s)|))ds=\int_{-\infty}^{\infty}f(x)L^{*}(x,t)dx.
\end{equation}
\end{proposition}
\begin{proof}
 By equations (\ref{stable-integral}) and
(\ref{local-definition})
\begin{eqnarray}
\int_{-\infty}^{\infty}f(x)L^{*}(x,t)dx & = &
\int_{-\infty}^{\infty}f(x)\int_{0}^{\infty}\bar{L}_{2}(s,t)d_{s}L_{1}(x,s)dx\nonumber\\
& =& \int_{0}^{\infty} \bar{L}_{2}(s,t)
d_{s}\int_{-\infty}^{\infty}f(x)L_{1}(x,s)dx\nonumber\\
& =& \int_{0}^{\infty}\bar{L}_{2}(s,t)d_{s}\int_{0}^{s}f(X(u))du\nonumber\\
& =& \int_{0}^{\infty} \bar{L}_{2}(s,t)f(X(s))ds\nonumber\\
&=& \int_{0}^{t}f(X(|Y(s)|))ds.
\end{eqnarray}
Hence we have the equation (\ref{iterated-int}). The joint
continuity of $ L^{*}(x,t)$ follows from the joint continuity of
the local times of Brownian motion and of symmetric stable
process.
\end{proof}
We now give the scaling property of local time of $Z^{1}$.
\begin{theorem}\label{scaling}
\begin{equation}\label{scaling-eq}
L^{*}(xt^{1/ 2 \alpha}, t)/
t^{1-1/2\alpha}\stk{(d)}{=}\int_{0}^{\infty}\bar{L}_{2}(s,1)d_{s}L_{1}(x,s)dx=L^{*}(x,1),\
\ x\in \RR{R}.
\end{equation}
\end{theorem}
\begin{corollary}
For each fixed $t\geq 0$, we have
\begin{equation}
L^{*}(0,t)/t^{1-1/2\alpha}\stk{(d)}{=}L^{*}(0,1).
\end{equation}
\end{corollary}

\begin{proof}[Proof of Theorem \ref{scaling}]
The following scaling properties of the Brownian local time and
stable local time are well-known.
\begin{equation}
\{L_{1}(x,t);\ x\in\RR{R},t\geq 0 \}\stk{(d)}{=}\{
\frac{1}{c^{1/2}}L_{1}(c^{1/2}x,ct); \ x\in\RR{R},t\geq 0\},
\end{equation}
and
\begin{equation}
\{L_{2}(x,t);\ x\in\RR{R},t\geq 0 \}\stk{(d)}{=}\{
\frac{1}{c^{1-1/\alpha}}L_{1}(c^{1/\alpha}x,ct); \
x\in\RR{R},t\geq 0\},
\end{equation}
where $c>0$ is an arbitrary fixed number. Consequently we have
\begin{eqnarray}
L^{*}(x,t)&=&\int_{0}^{\infty}\bar{L}_{2}(s,t)d_{s}L_{1}(x,s)\nonumber\\
&\stk{(d)}{=}&\frac{1}{c^{1-1/\alpha}}\int_{0}^{\infty}\bar{L}_{2}(c^{1/\alpha}s,ct)d_{s}L_{1}(x,s)\ \ c>0 \ \ \text{fixed} \nonumber\\
& \stk{(d)}{=} &
t^{1-1/\alpha}\int_{0}^{\infty}\bar{L}_{2}(\frac{s}{t^{1/\alpha}},1)d_{s}L_{1}(x,s),\
\ c=1/t, t>0 \ \ \text{fixed} \nonumber\\
 & \stk{(d)}{=} &
t^{1-1/\alpha}\int_{0}^{\infty}\bar{L}_{2}(u,1)d_{u}L_{1}(x,ut^{1/\alpha}),
\ \ u=s/t^{1/\alpha}, t>0 \ \ \text{fixed}\nonumber\\
 & \stk{(d)}{=} &
t^{1-1/\alpha}t^{1/2\alpha}\int_{0}^{\infty}\bar{L}_{2}(u,1)d_{u}L_{1}(\frac{x}{t^{1/\alpha}},u),\
\ t>0   \ \ \text{fixed}\ \ x\in \RR{R}.\nonumber
\end{eqnarray}
Clearly, the last equation is equivalent to (\ref{scaling-eq}).
\end{proof}

\subsubsection{Preliminaries}

In this section we will prove a lemma which is crucial in the
proof of the following theorems. The usual LIL or Kolmogorov's LIL
for Brownian motion which replaces the time parameter was used
essentially in the results in \cite{CCFR} and \cite{shi-yor} to
prove Kesten's LIL for iterated Brownian motion. However, there
does not exist an LIL of this type for symmetric $\alpha$-stable
process. To overcome this difficulty with the use of the following
lemma we show below that the LIL for the range process of
symmetric $\alpha$-stable process suffices to prove Theorem
\ref{local-liminf}.
\begin{lemma}\label{distr-measure}
Let $A\subset \RR{R}_{+}$ be Lebesgue measurable. Let $L(x,A)$ be
local time of Brownian motion over the set $A$. Then
$$
\sup_{x\in \RR{R}}L(x,A)\stk{(d)}{=}\sup_{x\in
\RR{R}}L(x,|A|)=\sup_{x\in \RR{R}}L(x,[0,|A|]),
$$
where $|.|$ denotes Lebesgue measure and  $\stk{(d)}{=}$ means
equality in distribution.
\end{lemma}
\begin{proof}
We use monotone class theorem from \cite{rev-yor}. Define
$$
\mathcal{S}=\{A\subset \RR{R}_{+}: \ \ \sup_{x\in
\RR{R}}L(x,A)\stk{(d)}{=}\sup_{x\in \RR{R}}L(x,|A|) \}.
$$
Obviously $\RR{R}_{+}\in \mathcal{S}$. Let $A,B\in \mathcal{S}$
and $A\subset B$. Since $L(x,A)$ is an additive measure in the set
variable
$$
\sup_{x\in \RR{R}}L(x,B)=\sup_{x\in \RR{R}}L(x,A) +\sup_{x\in
\RR{R}}L(x,B\setminus A),
$$
so
\begin{eqnarray}
\sup_{x\in \RR{R}}L(x,|B|) & \stk{(d)}{=}  & \sup_{x\in \RR{R}}L(x,B) \nonumber\\
&  = &\sup_{x\in \RR{R}}L(x,A)+\sup_{x\in \RR{R}}L(x,B\setminus
A).\nonumber
\end{eqnarray}
On the other hand,
$$
\sup_{x\in \RR{R}}L(x,|B|)=\sup_{x\in \RR{R}}L(x,|A|) +\sup_{x\in
\RR{R}}L(x,(|A|, |B|)),
$$
hence
\begin{eqnarray}
\sup_{x\in
\RR{R}}L(x,|B|) & \stk{(d)}{=} & \sup_{x\in \RR{R}}L(x,|A|)+\sup_{x\in \RR{R}}L(x,(|A|, |B|))\nonumber\\
& \stk{(d)}{=} & \sup_{x\in \RR{R}}L(x,A) +\sup_{x\in
\RR{R}}L(x,B\setminus A).\nonumber
\end{eqnarray}
For $0<|A|<\infty$, the moment generating function
of $\sup_{x\in \RR{R}}L(x,|A|)$ satisfies for some $\delta >0$,
(see \cite[Remark p.452]{kesten} )
$$0 <MGF(t)=E[e^{t\sup_{x\in \RR{R}}L(x,|A|)}]< \infty, \ \  for\ \ |t|<\delta . $$
Since $\sup_{x\in \RR{R}}L(x,|A|)$ and $\sup_{x\in
\RR{R}}L(x,(|A|, |B|))$ are independent and  similarly $\sup_{x\in
\RR{R}}L(x,A)$ and $\sup_{x\in \RR{R}}L(x,B\setminus A)$ are
independent, considering moment generating functions (in case $|A|=0$ or $|B\setminus
A|=0$, we do not need generating functions) which is the product of the moment generating functions, we get that
$B\setminus A \in \mathcal{S}$.

Let $A_{n}\subset A_{n+1}$ be an increasing sequence of sets in
$\mathcal{S}$. For $\lambda >0$,
\begin{eqnarray}
P[\sup_{x\in \RR{R}}L(x,\cup_{n=1}^{\infty}A_{n})\leq
\lambda]&=&P[\sup_{x\in \RR{R}}\lim_{n\to\infty}L(x,A_{n})\leq
\lambda]\nonumber\\
& = & \lim_{n\to\infty}P[\sup_{x\in \RR{R}}L(x,A_{n})\leq
\lambda]\nonumber\\
& = & \lim_{n\to\infty}P[\sup_{x\in \RR{R}}L(x,|A_{n}|)\leq
\lambda]\nonumber\\
& = & P[\sup_{x\in \RR{R}}L(x,|\cup_{n=1}^{\infty}A_{n}|)\leq
\lambda].\nonumber
\end{eqnarray}
Hence $ \cup_{n=1}^{\infty}A_{n} \in \mathcal{S}$.

Now to complete the proof we show that open intervals are in
$\mathcal{S}$. Every interval is in
$\mathcal{S}$, since the increments of Brownian motion are
stationary. It is clear that sets of measure zero are also in
$\mathcal{S}$. Hence $\mathcal{S}$ contains every Lebesgue
measurable set by monotone class theorem.

\end{proof}

\subsection{On upper-upper classes}

For further information on the L\'{e}vy classes we refer to
R\'{e}v\'{e}sz \cite{revez}.

In this section we prove
\begin{theorem}\label{uu}
There exists a $t_{0}=t_{0}(w)$ and a universal constant $C_{uu} \in (0,\infty )$
such that for $t>t_{0}$ we have
\begin{equation}\label{uu-bound}
L^{*}(0,t)\leq \sup_{x\in \RR{R}} L^{*}(x,t)\leq C_{uu}
t^{1-1/2\alpha}(\log \log t)^{1+1/2\alpha} \ \ a.s.
\end{equation}
\end{theorem}
\begin{proof}
By the LIL for the range $R(t)=|\{ x;\ Y(s)=x \ \ \text{for some
}\ \ s\leq t\}|$ given in Griffin \cite{griffin}

\begin{equation}\label{range-LIL-Y}
\limsup_{t\to\infty}t^{-1/\alpha}(\log \log
t)^{-(1-1/\alpha)}R(t)=c_{1}\ a.s.
\end{equation}
and the Kesten type LIL for $Y$, established by Donsker and
Varadhan \cite{don-var}, for some finite constant $c_{2}$

\begin{equation}\label{kesten-lil-stable}
\limsup_{t\to \infty}t^{-(1-1/\alpha)}(\log \log
t)^{-1/\alpha}\sup_{x\in \RR{R}} L_{2}(x,t)=c_{2} \ \ a.s.
\end{equation}
\begin{eqnarray}
 \sup_{x\in \RR{R}} L^{*}(x,t)&=& \sup_{x\in
\RR{R}}\int_{0}^{\infty}\bar{L}_{2}(s,t)d_{s}L_{1}(x,s)\nonumber\\
& = & \sup_{x\in \RR{R}}
\int_{Y[0,t]}\bar{L}_{2}(s,t)d_{s}L_{1}(x,s)\ \
a.s\nonumber\\
& = & O(t^{1-1/\alpha}(\log \log t)^{1/\alpha}\sup_{x\in
\RR{R}}L_{1}(x,\RR{R}_{+}\cap Y[0,t]))\nonumber\\
& = & O(t^{1-1/\alpha}(\log \log t)^{1/\alpha}\sup_{x\in
\RR{R}}L_{1}(x,|Y[0,t]|))\  (by\ Lemma \ \ref{distr-measure})\nonumber\\
& = &  O(t^{1-1/\alpha}(\log \log t)^{1/\alpha}\sup_{x\in
\RR{R}}L_{1}(x,ct^{1/\alpha}(\log \log t)^{1-1/\alpha}
)) \ \ a.s.\nonumber\\
&=& O(t^{1-1/2\alpha}(\log \log t)^{1+1/2\alpha}) \ \
a.s.\label{up-up-low}
\end{eqnarray}
with $t_{0}$ big enough, by using the LIL for the range of $Y$ and
LIL of Kesten  for $\sup_{s\in \RR{R}}\bar{L}_{2}(s,t)$ from
(\ref{range-LIL-Y})-(\ref{kesten-lil-stable}) respectively, and
then applying the Kesten's LIL once again to $\sup_{x\in
\RR{R}}L_{1}(x,ct^{1/\alpha}(\log \log t)^{1-1/\alpha} )$ given in
\cite{kesten}.
\end{proof}

\subsection{On upper-lower classes}
In this section we prove
\begin{theorem}
There exists a universal constant $C_{ul}\in (0,\infty)$ such that
\begin{equation}\label{ul-bound}
P[\sup_{x\in \RR{R}} L^{*}(x,t) \geq C_{ul} t^{1-1/2 \alpha} (\log
\log t)^{(1+\alpha)/2\alpha} \ \ i.o]=1.
\end{equation}

\end{theorem}
 Since the $\log \log $ powers do not match in equations
 (\ref{uu-bound}) and (\ref{ul-bound}), we cannot deduce an LIL for $ \sup_{x\in \RR{R}}
 L^{*}(x,t)$.
 \begin{proof}
 Proof follows from Theorem \ref{maintheorem1} and the
 observation
$$
t=\int_{x\in S(t)}L^{*}(x,t)dx\leq R^{*}(t)\sup_{x\in \RR{R}}
L^{*}(x,t), \ \ a.s.
$$
where $R^{*} (t)=|S(t) |=|\{ x:\ Z^{1}(s)=x \ \text{for some }\
s\leq t\}|$ and the fact that $R^{*}(t)\leq  2\sup_{0\leq s\leq
t}|Z^{1}_{s}|$.
\end{proof}

\subsection{On lower-upper classes}

In this section we prove

\begin{theorem}\label{lu}
There exists a universal constant $C_{lu} \in (0,\infty )$ such
that
\begin{equation}
P[L^{*}(0,t)\leq \sup_{x\in \RR{R}} L^{*}(x,t)\leq
C_{lu}(\frac{t}{\log \log t})^{1-1/2\alpha}\ \ \text{i.o}]=1.
\end{equation}
\end{theorem}
\begin{proof}
We have from Cs\'{a}ki and F\"{o}ldes \cite{cs-fold}: for $0\leq
a\leq 1$
\begin{equation}\label{small-BM}
P[\sup_{x\in \RR{R}} L_{1}(x,1)\leq a]\geq \exp(-\frac{c}{a^{2}}),
\end{equation}
for some absolute constant $c>0$. A similar result for the
local time of $Y$ is given in \cite{griffin}: there exists $\theta>0$ and
$c_{1}>0$ such that for $t$ large
\begin{equation}\label{small-stable}
P[\sup_{x\in \RR{R}} L_{2}(x,1)\leq \theta /(\log \log
t)^{1-1/\alpha}]\geq c_{1} \beta^{(\log \log t)},
\end{equation}
with $e^{-1} < \beta <1$.

 Define $c_{2}=\sqrt{4c/d}$  with $C_{\beta}=d+\log \beta^{-1} <1$ and
$t_{k}=\exp(k^{p})$, with $p\in (1, 1/C_{\beta})$ for $k$ large.
Consider
\begin{eqnarray}
s_{k}& \stk{def}{=} & 2c_{4}t_{k}^{1/\alpha}(\log \log
t_{k})^{1-1/\alpha},
\ \ c_{4} \ \ \text{constant in LIL of range of Y},\nonumber\\
D_{k} & \stk{def}{=} &  \left\{ \sup_{ x\in \RR{R} }
(L_{1}(x,s_{k})-L_{1}(x,s_{k-1}))\leq
\frac{c_{2}t_{k}^{1/2\alpha}}{2(\log \log
t_{k})^{1/2\alpha}} \right\},\nonumber\\
E_{k} & \stk{def}{=} &  \left\{ \sup_{ x\in \RR{R} }
(L_{2}(x,t_{k})-L_{2}(x,t_{k-1}))\leq \frac{\theta t_{k}^{1-
1/\alpha}}{(\log \log
t_{k})^{1-1/\alpha}} \right\},\nonumber\\
F_{k}& \stk{def}{=} & D_{k} \cap E_{k}.\nonumber
\end{eqnarray}
We have by means of (\ref{small-BM}) and (\ref{small-stable}),
$$
P[D_{k}]\geq P[\sup_{ x\in \RR{R} } L_{1}(x,1)\leq
\frac{c_{2}}{2(\log \log t_{k})^{1/2}} ]\geq \exp (-d \log \log
t_{k}),
$$
$$
P[E_{k}]\geq P[\sup_{ x\in \RR{R} } L_{2}(x,1)\leq
 \frac{\theta }{(\log \log t_{k})^{1-1/\alpha}}]\geq c_{1} \beta
^{(\log \log t_{k})}.
$$
Hence
$$
P[F_{k}]=P[D_{k}]P[E_{k}]\geq \frac{c_{1}}{k^{pC_{\beta}}},
$$
which implies $\sum_{k} P[F_{k}]=\infty$. Thanks to the
independence of the $F_{k}'s$, we can apply the Borel-Cantelli
lemma to conclude that, almost surely there exists infinitely many
$k$'s for which $F_{k}$ is realized. On the other hand, by the
Kesten LIL, for $X$ and $Y$ for all large $k$,
$$
\sup_{ x\in \RR{R} } L_{1}(x,s_{k-1})\leq 2 (s_{k-1} \log \log
s_{k-1})^{1/2}\leq \frac{c_{2}t_{k}^{1/2\alpha}}{(\log \log
t_{k})^{1/2\alpha}},
$$
$$
\sup_{ x\in \RR{R} } L_{2}(x,t_{k-1})\leq
2c_{5}t_{k-1}^{1-1/\alpha}(\log \log t_{k-1})^{1/\alpha}\leq
\frac{c_{6}t_{k}^{1-1/\alpha}}{(\log \log t_{k})^{1-1/\alpha}}.
$$
Therefore there exist infinitely many $k$'s such that
\begin{equation}\label{lu-up1}
\sup_{ x\in \RR{R} } L_{1}(x,s_{k})\leq
\frac{2c_{2}t_{k}^{1/2\alpha}}{(\log \log t_{k})^{1/2\alpha}},
\end{equation}
\begin{equation}\label{lu-up2}
\sup_{ x\in \RR{R} } L_{2}(x,t_{k})\leq
\frac{(\theta+c_{6})t_{k}^{1-1/\alpha}}{(\log \log
t_{k})^{1-1/\alpha}}.
\end{equation}
For those $k$ satisfying (\ref{lu-up1})-(\ref{lu-up2}), we have,
by the usual LIL for the range of $Y$ given in
(\ref{range-LIL-Y}),
\begin{eqnarray}
\sup_{ x\in \RR{R} } L^{*}(x,t_{k})& = & \sup_{ x\in \RR{R}
}\int_{0}^{\infty}\bar{L}_{2}(u,t_{k})d_{u}L_{1}(x,u)\nonumber\\
&=& \sup_{ x\in \RR{R}
}\int_{Y[0,t_{k}]}\bar{L}_{2}(u,t_{k})d_{u}L_{1}(x,u)\nonumber\\
& \leq  & 2 \sup_{ y\in \RR{R} }L_{2}(y,t_{k})\sup_{ x\in \RR{R}
}L_{1}(x,\RR{R}_{+} \cap Y[0,t_{k}])\nonumber\\
& \leq & 4c_{2}(\theta+c_{6})\left(\frac{t_{k}}{\log \log
t_{k}}\right)^{1-1/2\alpha}\ \ (by \ Lemma\
\ref{distr-measure}).\nonumber
\end{eqnarray}
\end{proof}

\subsection{The range}

In this section we will prove an LIL for the range $R^{*}(t)=| \{
x:\ Z^{1}(s)=x \ \text{for some }\  s\leq t\}|$. The idea of the
proof is to look at the large jumps of the symmetric stable
process which replaces the time parameter in the process
$Z^{1}(t).$ To prove LIL for the range of $Z^{1}$ we need several
lemmas. We adapt the arguments of Griffin \cite{griffin} to our
case in the following lemmas.

If $Y(t)$ is a process and $T$ is some, possibly random, time then
$Y^{*}(T)=\sup \{|Y(r)|: \ 0\leq r \leq T\}$. If $S< T$ then
$(Y(T)-Y(S))^{*}=\sup \{|Y(r)-Y(S)|: S\leq r\leq T\}.$
\begin{definition}
$T_{Y}(a)=\inf \{ s: \ |Y(s)-Y(s-)|>a \}$.
\end{definition}
We will usually write $T_{Y}(a)=T(a)$ if it is clear which process
we are referring to.

\begin{lemma}[Griffin \cite{griffin}]
The random variables $a^{-1/\alpha}Y^{*}(T(a^{1/\alpha}-))$ and
$Y^{*}(T(1)-)$ have the same distribution.
\end{lemma}
Using the scaling of Brownian motion it is easy to deduce
\begin{lemma}\label{Z-scale}
The random variables
$$a^{-1/2\alpha}X^{*}(Y^{*}(T(a^{1/\alpha}-)))$$ and
$$X^{*}(Y^{*}(T(1)-) )$$ have the same distribution.
\end{lemma}
 As in Griffin \cite{griffin}, we can decompose $Y$ as the sum of
 two independent L\'{e}vy processes
 $$
Y(t)=Y_{1}(t)+Y_{2}(t),
 $$
where
\begin{eqnarray}
Y_{2}& = &\sum_{s\leq t}(Y(s)-Y(s-))1\{
|Y(s)-Y(s-)|>1\}\nonumber\\
Y_{1}(t)&=& Y(t)-Y_{2}(t) .\nonumber
\end{eqnarray}
The L\'{e}vy measure of $X_{1}$ is given by $1\{|x|\leq
1\}|x|^{-1-\alpha}dx$ and the moment generating function by
$$
E[\exp(aY_{1}(t))]=\exp(t\psi (a)),
$$
where
$$
\psi (a)=\int_{-1}^{1}(e^{ax}-1)\frac{dx}{|x|^{1+\alpha}}.
$$
Observe that $\psi (a)\to 0$ as $a\to 0$.
\begin{lemma}\label{exp-stable}[Griffin \cite{griffin}]
If $a$ is small enough that $\psi (a)<2\alpha ^{-1}$, then
$$
E[\exp (a Y^{*}(T(1)-))]\leq \frac{8\alpha
^{-1}}{2\alpha^{-1}-\psi (a)}.
$$
\end{lemma}
We deduce the following from the last lemma.
\begin{lemma}\label{alpha-exp}
If $a$ is small enough that $\psi (a^{2}/2)\leq 2\alpha ^{-1}$,
then
$$
E[\exp ( a X^{*}(Y^{*}(T(1)-)))]\leq \frac{32\alpha
^{-1}}{2\alpha^{-1}-\psi (a^{2}/2)}.
$$
\end{lemma}
\begin{proof}
The moments of Brownian motion $X$ are given by $E[\exp (\theta
X(t))]=\exp  (\theta ^{2}t/2), $ so
$$
E[\exp ( a X^{*}(Y^{*}(T(1)-)))]=\int_{0}^{\infty}
\int_{0}^{\infty}E[ \exp (aX^{*}(l))]f_{l}(s) 2\alpha ^{-1}
e^{-2\alpha ^{-1}}dl ds,
$$
where $f_{l}(s)$ is the density of $Y^{*}(s)$.

 Now $P[X^{*}(t)>x]\leq 2 P[|X(t)|>x]$ for each $x>0$ $t>0$, hence
$$
E[ \exp (aX^{*}(l))]\leq 2 E[ \exp (a|X(l)|)]\leq  4 E[ \exp
(aX(l))]=\exp (a^{2}l/2),
$$
therefore
\begin{equation}\label{transition-1}
E[\exp ( a X^{*}(Y^{*}(T(1)-)))]\leq  4\int_{0}^{\infty}E[\exp
(\frac{a^{2}}{2}Y^{*}(s)) ]2\alpha ^{-1} e^{-2\alpha ^{-1}} ds.
\end{equation}
From Lemma \ref{exp-stable} we deduce that
 $$
E[\exp ( a X^{*}(Y^{*}(T(1)-)))]\leq \frac{32\alpha
^{-1}}{2\alpha^{-1}-\psi (a^{2}/2)}.
 $$
\end{proof}

\begin{definition}
$J(t,\gamma (t))= \# \{ s\leq t : \ \ |Y(s)-Y(s-)|>\gamma (t)\}$
where $\gamma (t)=(t/ \log \log t)^{1/\alpha}$.
\end{definition}

 We know from \cite{griffin} that  $ J(t,\gamma (t))$ has Poisson distribution with parameter $2\log \log t /\alpha$.
\begin{definition}
Fix $t>0$ and define
\begin{eqnarray}
T_{1}& = & \inf \{ s:\  |Y(s)-Y(s-)|\geq \gamma (t)\}\nonumber\\
T_{k+1}& =& \inf \{ s>T_{k}:\  |Y(s)-Y(s-)|\geq \gamma (t)\}.\nonumber\\
V_{k}^{1}& = & (X(|Y(T_{k}-)|)-X(|Y(T_{k-1})|))^{*}\nonumber \\
V_{k} &= & X^{*}(((Y(T_{k}-))-Y(T_{k-1}))^{*})\nonumber \\
W_{k}^{1} & = & (t/\log \log t)^{-1/2\alpha}V_{k}^{1}\nonumber\\
 W_{k} & = & (t/\log \log t)^{-1/2\alpha}V_{k}.\nonumber
\end{eqnarray}
\end{definition}
Observe that $W_{k}$, $k=1,2,\cdots$ are identically distributed
 as
$$X^{*}(Y^{*}(T(1)-))$$ by Lemma \ref{Z-scale}.

Observe also that $W_{k}^{1}$, $k=1,2,\cdots$ are identically
distributed as
$$
X^{*}((|Y(T_{2}-)|-|Y(T_{1})|)^{*}).
$$
Furthermore
$$
X^{*}((|Y(T_{2}-)|-|Y(T_{1})|)^{*})\leq
X^{*}((Y(T_{2}-)-Y(T_{1}))^{*}).
$$
Finally, observe that $ X^{*}((Y(T_{2}-)-Y(T_{1}))^{*})$ and $
X^{*}(Y^{*}(T(1)-))$ are identically distributed by Lemma
\ref{Z-scale}.

 Since the paths of $Y$ are not non-decreasing, the processes $W_{k}$, $k=1,2,\cdots$ are not independent.
 To get independent processes
  we have to disjointify the image of $Y$. So we define
 \begin{eqnarray}
 V^{*}_{1} & = &  \sup_{0\leq s\leq \sup_{0\leq r \leq T_{1}-}|Y(r)|} |X(s)|\nonumber\\
 V^{*}_{k} & = & \sup_{s,l\in A_{k}}|X(s)-X(l)|,\nonumber
 \end{eqnarray}
 where $A_{k}=|Y|[T_{k-1}, T_{k}-]\cap (|Y|[0, T_{k-1}-])^{C}$, $k=2,3 ,\cdots  $
 Observe that given $Y$, $V^{*}_{k}$ are independent for $k=1,2,\cdots  $, and $V^{*}_{k}\leq 2 V_{k}^{1}$. Define
  $$W_{k}^{*}  =  (t/\log \log t)^{-1/2\alpha}V_{k}^{*}.$$

 Now let $\varphi (t)$ denote the function  $(t/\log
\log t)^{1/2\alpha}\log \log t$.
\begin{lemma}\label{prob-bound}
If $\lambda$ is sufficiently large, then
$$
P[V_{1}^{*}+\cdots +V_{J(t,\gamma (t))+1}^{*}\geq \lambda \varphi
(t)]\leq (\log t)^{-2}.
$$
\end{lemma}

\begin{proof}
We have by a lemma in Griffin \cite{griffin}, for $\beta $
sufficiently large
\begin{equation}\label{number-of}
P[J(t,\gamma(t))\geq [\beta \log \log t]]\leq
\frac{1}{2}\frac{1}{( \log t )^{2}}.
\end{equation}
Since $V^{*}_{k}\geq 0$,
\begin{eqnarray}
P[V_{1}^{*}+\cdots +V_{J(t,\gamma (t))+1}^{*} \geq  \lambda
\varphi (t)]& \leq  &P[V_{1}^{*}+\cdots +V_{[\beta \log \log t]}^{*
}\geq \lambda \varphi
(t)]\nonumber\\
& \ \ + &  P[J(t,\gamma(t))\geq [\beta \log \log t]].\nonumber
\end{eqnarray}
Let $\xi =32\alpha ^{-1}/(2\alpha^{-1}-\psi ((2a)^{2}/2))$ and
choosing $\beta$ to satisfy (\ref{number-of}), we see that by
Lemma \ref{alpha-exp} for $a$ sufficiently small and  $ll t$
denoting $\log \log t$
\begin{eqnarray}
P[V_{1}^{*}+\cdots +V_{[\beta ll t]}^{*} &\geq  &\lambda \varphi
(t)] = P[W_{1}^{*}+\cdots +W_{[\beta ll t]}^{*}\geq \lambda
ll
t]\nonumber \\
& \leq & \exp(-a \lambda ll t)E[(E[\exp (aW_{1}^{*})|Y] \nonumber \\
 \ \  & & \times E[\exp (aW_{2}^{*})|Y] \cdots E[\exp (aW_{[\beta
 ll t]}^{*})|Y])]\nonumber\\
& \leq & \exp(-a \lambda ll t)E[(E[\exp (2aW_{1})|Y] \nonumber \\
 \ \  & & \times E[\exp (2aW_{2})|Y] \cdots E[\exp (2aW_{[\beta ll
t]})|Y])]\nonumber\\
& \leq & \exp(-a \lambda ll t)E[(4E[\exp (2a^{2}U_{1})|Y] \nonumber \\
 \ \  & & \times 4E[\exp (2a^{2}U_{2})|Y] \cdots 4E[\exp (2a^{2}U_{[\beta
 ll t]})|Y])]\nonumber\\
&\leq &\exp(-a \lambda ll t)(4E[\exp (2a^{2}U_{1})])^{[\beta ll
t]}\label{independent-Y}\\
&\leq &\exp(-ll t(a\lambda -\beta \log \xi))\nonumber\\
& \leq  & \frac{1}{2( \log t )^{2}},\nonumber
\end{eqnarray}
if $\lambda$ is sufficiently large. Where
$$U_{k}=(t/\log \log
t)^{-1/\alpha}(Y(T_{k}-)-Y(T_{k-1}))^{*}.$$  In the fourth line
inequality we use equation (\ref{transition-1}). In equation
(\ref{independent-Y}) we use the fact that $U_{k}'s$ are i.i.d.
with common distribution $Y^{*}(T(1)-)$ and Lemma
\ref{exp-stable}.
\end{proof}

\begin{theorem}\label{range-lower}
There exits a $t_{0}$ such that for $t>t_{0}$, and for  certain
constants $C,K \in (0,\infty )$
$$
 R^{*}(t)\geq C t^{1/2\alpha}(\log \log t)^{-(1+1/2\alpha )}\
 \text{a.s.}
$$
and
$$
P[R^{*}(t)\geq K t^{1/2\alpha}(\log \log t)^{1-1/2\alpha}\ \
\text{i.o}]=1.
$$
\end{theorem}

\begin{proof}
This follows from Theorem \ref{uu} and
$$
t=\int_{x\in S(t)}L^{*}(x,t)dx\leq R^{*}(t)\sup_{x\in \RR{R}}
L^{*}(x,t),
$$
where $R^{*}(t)=|S(t) |= |\{ x:\ Z^{1}(s)=x \ \text{for some }\
s\leq t\}|$. The last probability follows similarly using Theorem
\ref{lu}.
\end{proof}

\begin{proof}[Proof of Theorem \ref{range-upper}]
We will prove only the upper bound in the light of Theorem
\ref{range-lower}. To prove the upper bound observe that
$$
R^{*}(t)\leq V_{1}^{*}+\cdots +V_{J(t,\gamma (t))+1}^{*}.
$$
Thus by Lemma \ref{prob-bound} for $\lambda$ sufficiently large

$$
P[R^{*}(t)\geq \lambda \varphi (t)]\leq (\log t)^{-2}.
$$
Hence for large $n$,
\begin{eqnarray}
P[R^{*}(t)& \geq  &\lambda \varphi (t)\ \text{ for some } t\in
[2^{n},2^{n+1})]\nonumber\\
&\leq & P[R^{*}(2^{n+1})\geq \lambda \varphi(2^{n})]\nonumber\\
& \leq & P[R^{*}({2^{n+1}})\geq (\lambda c)   \varphi(2^{n+1})]\nonumber\\
& \leq  & \frac{C_{last}}{(n+1)^{2}}.\nonumber
\end{eqnarray}
if $\lambda$ is sufficiently large.  The result follows from Borel-Cantelli lemma.
\end{proof}

\subsection{On lower-lower classes}

In this section we prove

\begin{theorem}\label{local-lup}
There exists a $t_{0}=t_{0}(w)$  and a universal constant $C_{ll}\in (0,\infty )$ such that for $t>t_{0}$
\begin{equation}
\sup_{x\in \RR{R}}L^{*}(x,t)\geq C_{ll} (\frac{t}{\log \log
t})^{1-1/2\alpha}\ \ a.s.
\end{equation}
\end{theorem}
\begin{proof}
The proof follows from Theorem \ref{range-upper} and the
observation
$$
t=\int_{x\in S(t)}L^{*}(x,t)dx\leq R^{*}(t)\sup_{x\in \RR{R}}
L^{*}(x,t),
$$
where $R^{*}(t)=|S(t)| = |\{ x:\ Z^{1}(s)=x \ \text{for some }\
s\leq t\}|$.
\end{proof}
\begin{proof}[Proof of Theorem \ref{local-liminf}]
Theorems \ref{lu} and \ref{local-lup} imply the proof.
\end{proof}

 \textbf{Acknowledgments.} I would like to thank  Professor
Rodrigo Ba\~{n}uelos, my academic advisor, for  his guidance on
this paper.

\end{document}